\documentclass[aos,MSNbibl,nameyear,dvips]{arximspdf}
\usepackage{dcolumn}
\usepackage{mathrsfs}
\usepackage{graphicx}

\doi{10.1214/13-AOS1105} 
\volume{41}
\issue{2}
\pubyear{2013}
\firstpage{957}
\lastpage{981}

\makeatletter

\newcolumntype{d}[1]{D{.}{.}{#1}}

\newcommand{\lldots}{\cdots}

\newcommand{\cal}{\mathcal}

\newtheorem{theorem}{Theorem}[section]
\newtheorem{lemma}[theorem]{Lemma}

\newcommand{\one}{\mathbf{1}}
\newcommand{\Ef}{\mathbb E}
\newcommand{\Pf}{\mathbb P}

\makeatother

\begin{document}
\begin{frontmatter}

\title{Asymptotic theory with hierarchical autocorrelation:
Ornstein--Uhlenbeck tree models\thanksref{T1}}
\runtitle{Ornstein--Uhlenbeck hierarchical autocorrelation}

\thankstext{T1}{Supported in part by the NSF Grants
DEB-0830036 and DMS-11-06483.}

\begin{aug}
\author[A]{\fnms{Lam Si Tung} \snm{Ho}\ead[label=e1]{lamho@stat.wisc.edu}}
\and
\author[B]{\fnms{C\'ecile} \snm{An\'e}\corref{}\ead[label=e2]{ane@stat.wisc.edu}}
\runauthor{L. S. T. Ho and C. An\'e}
\affiliation{University of Wisconsin--Madison}
\address[A]{Department of Statistics\\
University of Wisconsin--Madison\\
1300 University Ave.\\
Madison, Wisconsin 53706\\
USA\\
\printead{e1}}
\address[B]{Department of Statistics\\
and\\
Department of Botany\\
University of Wisconsin--Madison\\
1300 University Ave.\\
Madison, Wisconsin 53706\\
USA\\
\printead{e2}} 
\end{aug}

\pdfauthor{Lam Si Tung Ho, Cecile Ane}

\received{\smonth{5} \syear{2012}}
\revised{\smonth{2} \syear{2013}}

%
\begin{abstract}
Hierarchical autocorrelation in the error term of linear models arises
when sampling units are related to each other according to a tree. The
residual covariance is parametrized using the tree-distance between
sampling units. When observations are modeled using an
Ornstein--Uhlenbeck (OU) process along the tree, the autocorrelation
between two tips decreases exponentially with their tree distance.
These models are most often applied in evolutionary biology, when tips
represent biological species and the OU process parameters represent
the strength and direction of natural selection. For these models, we
show that the mean is not microergodic: no estimator can ever be
consistent for this parameter and provide a lower bound for the
variance of its MLE. For covariance parameters, we give a general
sufficient condition ensuring microergodicity. This condition suggests
that some parameters may not be estimated at the same rate as others.
We show that, indeed, maximum likelihood estimators of the
autocorrelation parameter converge at a slower rate than that of
generally microergodic parameters. We showed this theoretically in a
symmetric tree asymptotic framework and through simulations on a large
real tree comprising 4507 mammal species.
\end{abstract}

%

\end{frontmatter}

\section{Introduction and overview of main results}

\subsection{Motivation}
\begin{figure}

\includegraphics{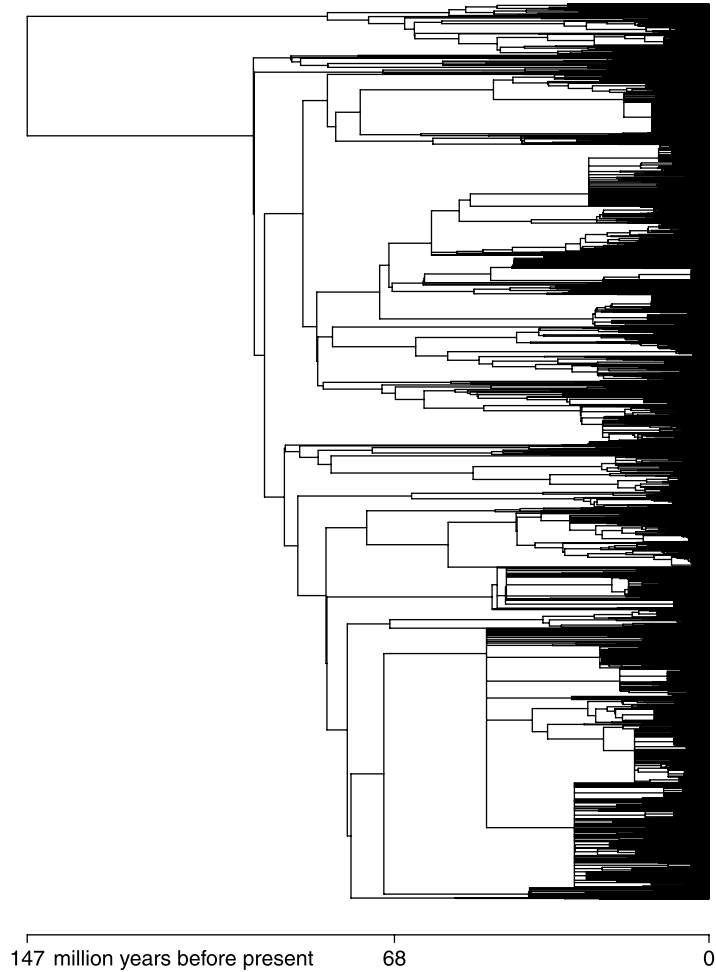}

\caption{Family tree of 4507 mammal species
[Bininda-Emonds
et~al. (\citeyear{binindaEmonds-etal07})]. Branch lengths indicate estimated
diversification times on the horizontal axis. The Cretaceous/Tertiary
mass extinction event marked the extinction
of dinosaurs 65.5 million years ago. 
Cooper and Purvis (\citeyear{cooperPurvis10}) used body mass data
available for 77\% of these species to infer the mode of evolution:
neutral evolution (BM) versus natural selection (OU).}
\label{figmammaltree}
\end{figure}

This work is motivated by the availability of very large data sets to
compare biological species, and by the current lack of asymptotic
theory for the models that are used to draw inference from species comparisons.
For instance, \citet{cooperPurvis10} studied the evolution of body size
in mammals using data from 3473 species
whose genealogical relationships are depicted by their family tree in
Figure~\ref{figmammaltree}. 
Even from this abundance of data, Cooper and Purvis found a lack of power
to discriminate between a model of neutral evolution versus a model with
natural selection. 
%
%
To model neutral evolution, body size is assumed to follow
a Brownian motion (BM) along the branches of the tree, with observations
made on present-day species at the tips of the tree.
To model natural selection, body size is assumed to follow an
Ornstein--Uhlenbeck (OU) process, whose parameters represent a
selective body size ($\mu$) and a selection strength ($\alpha$).
The lack of power observed by Cooper and Purvis suggests a nonstandard
asymptotic behavior of the model parameters, which is the motivation for
our work.

%
%
\begin{figure}[b]

\includegraphics{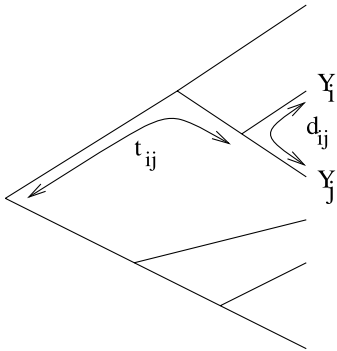}

\caption{Correlation (or residual correlation) between
observations at tips $i$ and $j$ are parametrized in the OU
model as a function of the tree distance $d_{ij}$ between $i$ and
$j$ and of the length $t_{ij}$ of their shared path from the root.
For instance, Cooper and Purvis (\citeyear{cooperPurvis10}) considered
body mass ($Y$) across 3473 mammal species ($i,j=1,\ldots,3473$).}
\label{figfig1}
\end{figure}
%

\subsection{Tree structured autocorrelation}
Hierarchical autocorrelation, as depicted in the mammalian tree,
arises whenever sampling units are related to each other through
a vertical inheritance pattern, like biological species,
genes in a gene family or human cultures.
In the genealogical tree describing the relatedness between
units, internal nodes represent ancestral unobserved units
(like species or human languages).
Branch lengths measure evolutionary time between branching events
and define a distance between pairs of sampling units.
This tree and its branch lengths can be used to parametrize the
expected autocorrelation. For doing so, the BM and the OU process
are the two most commonly used models.
They are defined as usual along each edge in the tree.
At each internal node, descendant lineages
inherit the value from the parent edge just prior to the
branching event, thus ensuring continuity of the process.
Conditional of their starting value, each lineage then evolves
independently of the sister lineages.
BM evolution of the response variable (or of error term)
along the tree results in normally distributed errors and in a
covariance matrix governed by the tree, its branch lengths
and a single parameter $\sigma^2$.
The covariance between two tips $i$ and $j$ is simply $\sigma^2 t_{ij}$,
where $t_{ij}$ is the shared time from the root of the tree to the tips
(Figure~\ref{figfig1}).
Under the more complex OU process,
changes toward a value $\mu$ are favored over
changes away from this value, making the OU model
appropriate to address biological questions about the
presence or strength of natural selection.
This model is defined by the following stochastic equation
[\citet{ikedaWatanabe81}]:
$dY_t = -\alpha(Y_t-\mu)\,dt + \sigma\,dB_t$
where $Y$ is the response variable (such as body size),
$\alpha$ is the selection strength and $B_t$ is a BM process.
In what follows, $\mu$ is called the ``mean'' even though
it is not necessarily the expectation of the observations.
It is the mean of the stationary distribution of the OU process,
and it is the mean at the tips of the tree if the state at the root
has mean $\mu$.
In the biology literature, $\mu$ is called the ``optimal'' value or
``adaptive optimum'' in reference to the action of natural selection,
but this terminology could cause confusion here with likelihood optimization.
The parameter $\alpha$ measures the strength of the pull back to $\mu$.
High $\alpha$ values result in a process narrowly distributed
around $\mu$, as expected under strong natural selection if
the selective fitness of the trait is maximized at
$\mu$ and drops sharply away from $\mu$.
Simple mathematical models of natural selection at the level of
individuals result in the OU process for the population mean
[\citet{lande79,hansenMartins96}].
If $\alpha=0$, the OU process reduces to a BM with
no pull toward any $\mu$ value, as if the trait under consideration
does not affect fitness.
While some applications focus on
the presence of natural selection $(\alpha\neq0)$ such as
\citet{cooperPurvis10}, other applications are
interested in models where $\mu$ takes
different values $(\mu_1,\ldots,\mu_p)$ along different
branches in the tree, to model different adaptation regimes
[e.g., \citet{butlerKing04}].
%
Other applications assume a randomly varying 
$\mu$ along the tree, varying linearly with explanatory variables
[\citet{hansenPienaarOrzack08}].
In our work, we develop an asymptotic theory for the simple
case of a constant $\mu$ over the whole tree.
The covariance between two observed tips
depends on how the unobserved response at the root is treated.
It is reasonable to assume that this value $y_0$ at the root is a
random variable with the stationary Gaussian distribution
with mean $\mu$ and variance $\gamma=\sigma^2/(2\alpha)$.
With this assumption, the observed process ${(Y_i)}_{i\in\mathrm
{tips}}$ is
Gaussian with mean $\mu$ and variance matrix
%
%
\begin{equation}
\label{eqrandomrootV}
\gamma\mathbf{V}\qquad\mbox{with }V_{ij}=e^{-\alpha d_{ij}},
\end{equation}
where $d_{ij}$ is the tree distance between tips $i$ and $j$,
that is, the length of the path between $i$ and $j$.
Therefore, the strength $\alpha$ of natural selection
provides a direct measure of the level of autocorrelation.
If instead we condition on the response value $y_0$ at the root,
the Gaussian process has mean
$(1-e^{-\alpha t_{ii}})\mu+ e^{-\alpha t_{ii}}y_0$ for tip $i$
and variance matrix
%
%
\begin{equation}
\label{eqfixedrootV} \gamma\mathbf{V}\qquad\mbox{with }
V_{ij}=e^{-\alpha d_{ij}}
\bigl(1-e^{-2\alpha t_{ij}}\bigr),
\end{equation}
where, again, $t_{ii}$ is the distance from the root to tip $i$, and
$t_{ij}$ is the shared time from the root to tips $i$ and $j$
(Figure~\ref{figfig1}).

\subsection{Main results and link to spatial infill asymptotics}
In contrast to autocorrelation in spatial data or time series,
hierarchical autocorrelation has been little considered in the
statistics literature, even though tree models have been
used in empirical studies for over 25 years.
The usual asymptotic properties have mostly been taken for granted.
Recently, \citet{ane08} showed that the maximum likelihood (ML) estimator
of location parameters is not consistent under the BM tree model
as the sample size
grows indefinitely, proving that the basic
consistency property should not be taken for granted.
However, \citet{ane08} did not consider the more complex OU model,
for which the ML estimator admits no analytical formula.

In the spatial infill asymptotic framework when
data are collected on a denser and denser set of locations within a
fixed domain, $\sigma^2$ can be consistently estimated, but $\alpha$ cannot
under an OU spatial autocorrelation model in dimension $d \leq3$
[\citet{zhang04}].
Recently, $\alpha$ has been proved to be consistently estimated under
OU model when $d \geq5$ [\citet{anderes2010consistent}].
We uncover here a similar asymptotic behavior under the OU tree model.
Just like in infill asymptotics, the tree structure implies that
all sampling units may remain within a bounded distance of each
other, and that the minimum correlation between any pair of
observations does not go down to zero with indefinitely large
sample sizes. It is therefore not surprising that some properties
may be shared between these two autocorrelation frameworks.
Under infill asymptotics, microergodic parameters can usually be
consistently estimated [see \citet{zhang2005towards}] while nonmicroergodic
parameters cannot (e.g., $\alpha$). 
A parameter is microergodic when two different values for it
lead to orthogonal distributions for the complete,
asymptotic process [\citet{stein99}].

In Section~\ref{secmic}, we prove that the mean $\mu$ is
nonmicroergodic under the OU autocorrelation framework, and we provide
a lower bound for the variance of the MLE of $\mu$. We also give a
sufficient condition for the microergodicity of the OU covariance
parameters $\alpha$ and $\sigma^2$ (or $\gamma$) based on the
distribution of internal node ages.
The microergodic covariance parameter under spatial infill asymptotics
with OU autocorrelation, $\sigma^2$, is recovered as microergodic
if $0$ is a limit point of the sequence of node ages,
that is, with dense sampling near the tips.
Our condition for microergodicity suggests that some parameters
may not be estimated at the same rate as others. In Section~\ref{secrates},
we illustrate this theoretically for a symmetric tree asymptotic framework,
where we show that the REML estimator of $\alpha$ converges at a
slower rate than that of the generally microergodic parameter.
We also illustrate that the ML estimate convergence rate of $\alpha$
is slower than that of $\sigma^2$, through simulations
on a large 4507-species real tree showing dense sampling near the tips.

In most of this work, we only consider ultrametric trees, that is,
trees in which the root is at equal distance from all the tips.
This assumption is very natural for real data.
We also focus on model (\ref{eqrandomrootV}), because the model
matrix is not of full rank under model (\ref{eqfixedrootV}) on
an ultrametric tree.

\subsection{Other tree models in spatial statistics}
Trees have already been used for various purposes in spatial
statistics. When considering different resolution scales,
the nesting of small spatial regions into larger regions can be
represented by a tree. The data at a coarse scale for a given
region is the average of the observations at a finer scale
within this region. For instance,
\citet{huangCressieGabrosek02} use this ``resolution'' tree
structure to obtain consistent estimates at different scales,
and otherwise use a traditional spatial correlation
structure between locations at the finest level.
In contrast, the tree structure in our model is the fundamental
tool to model the correlation between sampling units, with
no constraint between values at different levels.
Trees have also been used to capture the correlation among
locations along a river network
[\citet{cressie-etal06,verhoef-etal06,verhoefPeterson10}, and discussion].
A~river network can be represented by a tree with the
associated tree distance.
To ensure that the covariance matrix is positive definite,
moving average processes have been introduced, either averaging
over upstream locations or over downstream locations, or both.
There are two major differences between our model and these
river network models.
First, the correlation among moving averages considered in
\citet{cressie-etal06} and \citet{verhoefPeterson10} decreases
much faster than the correlation considered in this work.
Most importantly, any location along the river is observable,
while observations can only be made at the leaves of
the tree in our framework.

\section{Microergodicity under hierarchical autocorrelation} \label{secmic}
The concept of microergodicity was formalized by \citet{stein99}
in the context of spatial models. This concept was especially
needed in the infill asymptotic framework, when some parameters
cannot be consistently estimated even if the whole process is observed.
Specifically, consider the complete process ${(Y_s)}_{s\in S}$
where $S$ is the space of all possible observation units. In spatial
infill asymptotics, $S$ can be the unit cube $[0,1]^d$.
In our hierarchical framework, we consider a sequence of
nested trees converging to a limit tree, which is the union of all
nodes and edges of the nested trees.
In this case, $S$ is the set of all tips in the limit tree.
Consider a probability model
${(P_\theta)}_{\theta\in\Theta}$ on ${(Y_s)}_{s\in S}$.
A function $f(\theta)$ of the parameter vector is said to be microergodic
if for all $\theta_1, \theta_2 \in\Theta$,
$f(\theta_1)\neq f(\theta_2)$ implies that
$P_{\theta_1}$ and $P_{\theta_2}$ are orthogonal.
If a parameter is not microergodic, then there is no hope of
constructing any consistent estimator for it; see \citet{zhang04} for
an excellent explanation.
In spatial infill asymptotics with OU correlation
in dimension $d \leq3$, $\alpha$ and $\gamma$ are not microergodic
even though $\sigma^2$ is [\citet{zhang04}], and the MLE of $\sigma
^2$ is strongly
consistent [\citet{ying91}].
Also note that the microergodicity of
$(\gamma,\alpha)$ is equivalent to the microergodicity of both
$\gamma$ and $\alpha$.

\subsection{Theory of equivalent Gaussian measures}
We recall here the theory of equivalent Gaussian measures,
which we apply to\vadjust{\goodbreak} Ornstein--Uhlenbeck tree models in the next section.
We consider two Gaussian measures $P_k$ $(k=1,2)$ on the $\sigma
$-algebra $\mathscr{U}$ generated by a sequence of random variables
$(Y_j)_{j=1}^\infty$, a linearly independent basis for both ${\cal
H}_1$ and ${\cal H}_2$ where ${\cal H}_k$ is the Hilbert space
generated by $(Y_j)_{j=1}^\infty$ with linear product: $\langle
Y_{j_1},Y_{j_2}\rangle= \operatorname{cov}_k(Y_{j_1} Y_{j_2})$
for $k=1$ or $2$. The entropy distance between equivalent Gaussian
measures $P_1$ and $P_2$ on the $\sigma$-algebra $\mathscr{U}'
\subset\mathscr{U}$ is defined as twice the symmetrized
Kullback--Leibler divergence,
\[
r\bigl(\mathscr{U}'\bigr) = - \biggl[ \Ef_{P_1} \log{
\frac{P_2(dw)}{P_1(dw)}} + \Ef_{P_2} \log{\frac{P_1(dw)}{P_2(dw)}}
\biggr].
\]
We will use the following properties proved in \citet{ibragimovRozanov78}:
%
%
\begin{equation}
\label{eqaddmic01} r\bigl(\mathscr{U}'\bigr) \leq r\bigl(
\mathscr{U}''\bigr) \qquad\mbox{for } \mathscr
{U}' \subset\mathscr{U}''.
\end{equation}
Consider nonsingular Gaussian measures $P_1$ and $P_2$ on the $\sigma
$-algebra $\mathscr{U}_n$ generated by $(Y_j)_{j=1}^n$. Let $r_n =
r(\mathscr{U}_n)$. Then $(r_n)_{n=1}^\infty$ is nondecreasing and
%
%
\begin{equation}
\label{thmaddmic01} P_1 \,\bot\, P_2
\quad\Leftrightarrow\quad r_n \to\infty\quad\mbox{and}\quad P_1 \equiv
P_2 \quad\Leftrightarrow\quad r_n \to r < \infty.
\end{equation}
We now recall how to calculate $r_n$ as described in \citet{stein99};
see also \citet{ibragimovRozanov78}. Consider a new basis
$(Y_{1,n},\ldots, Y_{n,n})$ obtained by linearly transforming
$(Y_1,\ldots, Y_n)$ such that this new basis is centered orthonormal
under $P_1\dvtx\Ef_1 Y_{j,n} =0$ and
$\operatorname{cov}_1(Y_{j_1,n},Y_{j_2,n}) = \delta_{j_1,j_2}$ is $1$
if $j_1=j_2$ and is $0$ otherwise, and such that
$\operatorname{cov}_2(Y_{j_1,n},Y_{j_2,n}) = \sigma^2_{j_1,n} \delta_{j_1,j_2}$
for some $\sigma^2_{j_1,n}$. Also set $m_{j,n} = \Ef_2 Y_{j,n}$.
Then
\[
r_n = \frac{1}{2} \sum_{j=1}^{n}
\bigl( \sigma^2_{j,n} +1/\sigma^2_{j,n}
-2 +m^2_{j,n} + m^2_{j,n}/
\sigma^2_{j,n} \bigr).
\]
\citet{rao1963discrimination} take a similar approach using the
Hellinger distance instead of the entropy distance $r_n$. They show
that the following condition
is sufficient for the orthogonality of $P_1$ and $P_2$:
%
%
\begin{equation}
\label{condRao} \lim_{n \to\infty} {\sum_{j=1}^{n}{
\bigl( \sigma^2_{j,n} - 1\bigr)^2}} = \infty.
\end{equation}

\subsection{Microergodicity of Ornstein--Uhlenbeck tree models}

We say that $\mathbb{T}$ is a subtree of tree $\mathbb{T}'$ if we can
get $\mathbb{T}$ by removing some branches from $\mathbb{T}'$. We
consider a nested sequence of trees $(\mathbb{T}_n)_{n=1}^\infty$
such that $\mathbb{T}_{n-1}$ is a subtree of $\mathbb{T}_n$ for every $n$.
This is\vspace*{1pt} to ensure that the observations $(Y_j)_{j=1}^n$ at the tips of
$\mathbb{T}_n$ provide a well-defined infinite sequence
${(Y_n)}_{n\geq1}$.
One essential assumption is that trees are ultrametric, that is, the
distance from the root to leaf nodes of tree $\mathbb{T}_n$ is assumed
to be the same for all tips. This is equivalent to saying that the tree
distances between tips define an ultrametric metric. This assumption
comes in naturally.\vadjust{\goodbreak} If the distance from the root to all tips is
constant, models (\ref{eqrandomrootV})
and (\ref{eqfixedrootV}) predict equal variances and equal means
at the tips, which are reasonable assumptions. Ultrametric trees
arise in most applications when tips are extant species sampled
at the present time, and branch lengths represent time calibrated
in millions of years, for instance. Define $\mathscr{I}^{\mathbb{T}_n}$
as the set of all internal nodes of tree $\mathbb{T}_n$ (including the root)
and $\mathscr{I} = \bigcup_{n=1}^\infty{\mathscr{I}^{\mathbb{T}_n}}$.
Let $(T_i)_{i \in\mathscr{I}}$ be the sequence of node ages. The age
of a node
is the distance from the node to any of its descendant tip.
This is well defined on ultrametric trees. $\mathscr{I}^{\mathbb{T}_n}$
is a subset of $\mathscr{I}^{\mathbb{T}_{n+1}}$ so $(T_i)_{i \in
\mathscr{I}}$ is a well-defined infinite sequence. In most of what
follows, we will assume that:
\begin{longlist}[(C)]
\item[(C)] $(\mathbb{T}_n)_{n=1}^\infty$ is a nested sequence of
ultrametric trees and the sequence of internal node ages $(T_i)_{i \in
\mathscr{I}}$ is bounded.
\end{longlist}
Without loss of generality, we can assume
that all trees are bifurcating because a multifurcating tree can be
made into
a bifurcating tree with some zero branch lengths.
With this assumption $\mathscr{I}^{\mathbb{T}_n}$ contains $n-1$
internal nodes.
This is equivalent to counting nodes and their ages with multiplicity,
where an internal node having $d$ descendants contributes
his age $d-1$ times.

Theorems~\ref{thmaddmica} and~\ref{thmaddmic} below state general
results on the microergodicity of parameters in OU tree models.
Our main tool is the equivalence (\ref{thmaddmic01}) applied to
$r_n=r(\mathbb{T}_n)$,
the entropy distance between $P_{\theta_1}$ and $P_{\theta_2}$ for
two parameter sets
$\theta_k=(\mu_k,\alpha_k,\gamma_k)$, $k=1,2$, on the $\sigma
$-algebra generated by $(Y_j)_{j=1}^n$.

\subsection{\texorpdfstring{Microergodicity of the mean $\mu$}{Microergodicity of the mean mu}}

%
\begin{theorem}\label{thmaddmica}
Under OU model (\ref{eqrandomrootV}) and condition \textup{(C)}, $\mu$ is not
microergodic.
\end{theorem}

The theorem follows directly from (\ref{thmaddmic01})
and the boundedness of $(T_i)_{i \in\mathscr{I}}$ once the following
upper bound is established:
%
%
\begin{equation}
\label{lemaddmic04} r(\mathbb{T}) \leq(\mu_1 -
\mu_2)^2/\bigl(\gamma_1 e^{- 2 \alpha_1 T}
\bigr),
\end{equation}
if $\alpha_1 = \alpha_2$ and $\gamma_1 = \gamma_2$,
where $T$ is the age of the root of $\mathbb{T}$
(Appendix~\ref{pfsecmiclemmas}).
%
%
%
\begin{figure}

\includegraphics{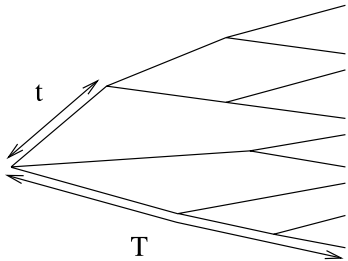}

\caption{Ultrametric tree with all tips at equal
distance $T$ from the root. The root has $k=3$ children here, and $t$
is the minimum distance from the root to its children.}
\label{figfig2}
\end{figure}
One consequence is that there is no consistent
estimator for $\mu$. To illustrate this, we consider the MLE of $\mu$
and provide a lower bound for its variance. We let $t$ be the length of
the shortest branch stemming from the root and $k$ the
number of daughters of the root (Figure~\ref{figfig2}).

%
\begin{theorem}\label{thminc01}
Assume OU model (\ref{eqrandomrootV}) on an ultrametric tree.
Let $\hat{\mu}$ be the MLE of $\mu$ conditional on
some possibly wrong value $\alpha_*$ of $\alpha$. Then
%
%
\begin{equation}
\label{eqninc01} \operatorname{var}(\hat{\mu}) \geq\frac{\sigma^2}{2
\alpha}
e^{-2
\alpha T} \biggl( 1 + \frac{e^{2 \alpha t} - 1}{k} \biggr).
\end{equation}
The equality holds if and only if $\alpha$ is known
($\alpha_*=\alpha$) and the tree is a star tree with the root as
unique internal node, in which case $k=n$ and $t=T$.
If $T$ is bounded as the sample size $n$ grows and $\alpha>0$,
then $\hat\mu$ is not consistent.
\end{theorem}

The second part of the theorem follows directly from the lower bound~(\ref{eqninc01}).
Note that $\hat{\mu}$ is Gaussian with mean $\mu
$. Therefore, the lower bound of its variance implies that $\hat{\mu
}$ cannot converge to $\mu$. Hence, it is not consistent.

The assumption that $\alpha>0$ is trivial. When $\alpha=0$, the
OU process reduces to a BM where $\mu$
has no influence on the process. In that case, $\mu$ is no longer
a parameter in the model.
As expected, the lower bound on the variance of $\hat\mu$ is heavily
influenced by the actual value of the correlation parameter $\alpha$.
The precision of $\hat\mu$ is weakest when autocorrelation is strong,
that is, when $\alpha$ is small, for a given value of $\gamma=\sigma
^2/(2\alpha)$.

\begin{figure}[b]

\includegraphics{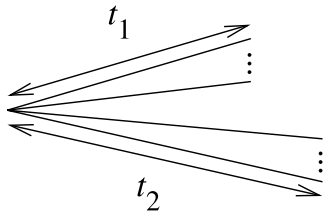}

\caption{Example of a nonultrametric tree on which
$\mu$ can be consistently estimated.}
\label{fignotultrametric}
\end{figure}

The ultrametric assumption is necessary. If the tree is not ultrametric,
model (\ref{eqfixedrootV})
predicts unequal variances and most importantly unequal means
at the tips. Such trees can carry more information about $\mu$.
Consider, for instance, the star tree in Figure~\ref{fignotultrametric},
in which all tips are directly connected to the root, by a branch of
%
%
length $t_1$ for half of the tips and of length $t_2$ for the other half
of the tips.
If $t_1\neq t_2$ the variance of $\hat\mu$ goes to $0$ as the sample size
grows (see Appendix~\ref{pfsecmiclemmas}), thus providing
a counterexample to Theorem~\ref{thminc01} when the ultrametric assumption
is violated.

\begin{pf*}{Proof of Theorem~\ref{thminc01}}
To prove (\ref{eqninc01}), we note that
$\hat{\mu} = (\one^{t} V_{\alpha_*}^{-1} \one)^{-1}\*
\one^{t} V_{\alpha_*}^{-1} Y$,
where $\one$ is a vector of ones.
This estimator is unbiased and has variance
$\frac{\sigma^2}{2 \alpha}(\one^{t} V_{\alpha}^{-1} \one)^{-1}$
when $\alpha_*=\alpha$ 
is known. Its variance is larger when $\alpha$ 
is unknown, by the Gauss--Markov theorem. For this reason, we only
need to prove the following lemma (which is done in Appendix \ref
{pfsecmiclemmas}).
\end{pf*}

%
\begin{lemma}\label{leminc02}
For all $\alpha>0$,
$
(\one^{t} V_{\alpha}^{-1} \one)^{-1} \geq e^{-2 \alpha T}
+ \frac{1}{k}(e^{-2 \alpha(T-t)} - e^{-2 \alpha T})
$
with equality if the tree is a star with $k$ branches stemming
from the root.
\end{lemma}

Theorem~\ref{thminc01} can be applied to any tree growth asymptotic framework,
so long as $T$ is bounded. For instance, both conditions are met almost
surely with $k=2$
under the coalescent model [Kingman (\citeyear{kingman82a,kingman82b})].
Even if these conditions do not hold asymptotically, 
(\ref{eqninc01}) provides a finite-sample upper bound on the
estimator's precision. This inequality can be used, for instance,
under the Yule model of tree growth [\citet{yule25,aldous01}]
if we let both $T$ and $n$ increase indefinitely.


\subsection{\texorpdfstring{Microergodicity of the autocorrelation parameter $(\gamma,\alpha)$}
{Microergodicity of the autocorrelation parameter (gamma, alpha)}}

%
\begin{theorem}\label{thmaddmic}
Under OU model (\ref{eqrandomrootV}) and condition \textup{(C)}:
%
\renewcommand\thelonglist{(\alph{longlist})}
\renewcommand\labellonglist{\thelonglist}
\begin{longlist}
\item\label{thmaddmicb} Let $t_0$ be a limit point of $(T_i)_{i \in
\mathscr{I}}$. Then
$f_{t_0}(\gamma,\alpha)$ is microergodic, where
\[
f_t(\gamma,\alpha) = \cases{\gamma\bigl(1-e^{-2 \alpha t}\bigr),
&\quad$t > 0$,
\cr
\gamma\alpha, &\quad$t=0$.}
\]
\item\label{thmaddmicc} If $\sum_{i \in\mathscr{I}}{(T_i - t)^2} =
\infty$ for all
$t\geq0$, then $(\gamma,\alpha)$ is microergodic.
Note that this condition is satisfied if $(T_i)_{i \in\mathscr{I}}$
has 2 or more
limit points.
\end{longlist}
\end{theorem}

\begin{pf}
The key idea is to reduce the tree for a lower bound of $r(\mathbb{T}_n)$.
We will consider subtrees that provide independent contrasts, sufficient
to ensure microergodicity.
Our constructive proof could be used to construct estimators
based on a restricted set of contrasts, but we do not pursue this here.
Let $i \in\mathscr{I}$ be an arbitrary internal node, and $Y^i_1$,
$Y^i_2$ be two leaves having $i$ as their most recent common ancestor.
Let $p_i$ be the path connecting $Y^i_1$ and $Y^i_2$. We define
$C^{p_i}_i = Y^i_1 - Y^i_2$ as a contrast with respect to internal node
$i$ and path $p_i$. For convenience, we define $T_{C^{p_i}_i} = T_i$.
The following lemma is proved in Appendix~\ref{pfsecmiclemmas}.

%
\begin{lemma}\label{lemaddmic06} We have that
$C^{p_i}_i \sim N(0,2\gamma(1-e^{-2\alpha T_i}))$.
Also, $C^{p_{i_1}}_{i_1}$ and $C^{p_{i_2}}_{i_2}$ are independent if
their paths $p_{i_1}$ and $p_{i_2}$ do not intersect.
\end{lemma}

\begin{pf*}{Proof of part~\ref{thmaddmicb}}
We denote
$\mathscr{I}^{\mathbb{T}}_S = \{ i\dvtx T_i \in S, i \in\mathscr
{I}^{\mathbb{T}} \}$
the set of internal nodes of $\mathbb{T}$ whose ages lie in $S$.
Let $(\gamma_1, \alpha_1)$ and $(\gamma_2,\alpha_2)$ such that
$f_{t_0}(\gamma_1,\alpha_1) \ne f_{t_0}(\gamma_2,\alpha_2)$. Denote
\[
g(t) = \frac{1}{2} \biggl(\frac{f_t(\gamma_1,\alpha_1)}{f_t(\gamma
_2,\alpha_2)} + \frac{f_t(\gamma_2,\alpha_2)}{f_t(\gamma_1,\alpha
_1)}-2 \biggr),\qquad
t \in\bigl[0,T^*\bigr],
\]
and let $\delta= g(t_0)/2 > 0$. Note that $g$ is continuous at $t_0$,
so there exists $\varepsilon_\delta> 0$ such that $g(t) \geq g(t_0) -
\delta$ for all $t$ satisfying $|t-t_0| < \varepsilon_\delta$.
We now use Lemma~\ref{lemaddmic05} (in Appendix~\ref{pfsecmicic})
to select a large set ${\mathscr C}_n$ of independent contrasts with
respect to internal nodes whose ages are in $(t_0 - \varepsilon_\delta,
t_0 + \varepsilon_\delta)$ such that $|{\mathscr C}_n| \geq\frac{1}{2}
|\mathscr{I}^{\mathbb{T}_n}_{(t_0 - \varepsilon_\delta, t_0 +
\varepsilon
_\delta)}|$. Let $r({\mathscr C}_n)$ be the entropy distance between
$P_{\theta_1}$ and $P_{\theta_2}$ on the $\sigma$-algebra generated
by ${\mathscr C}_n$. By (\ref{eqaddmic01}) and direct calculation,
\[
r(\mathbb{T}_n) \geq r(\mathscr{C}_n) = \sum
_{C \in\mathscr{C}_n} g(T_{C}) \geq|\mathscr{C}_n|
\bigl( g(t_0) - \delta\bigr) = \delta|\mathscr{C}_n|
\geq\frac{\delta}{2} \bigl|\mathscr{I}^{\mathbb{T}_n}_{(t_0-\varepsilon
_\delta,t_0+\varepsilon_\delta)}\bigr|.
\]
Clearly
$|\mathscr{I}^{\mathbb{T}_n}_{(t_0-\varepsilon_\delta,t_0+\varepsilon
_\delta)}|
\to\infty$ if $t_0$ is a limit point of $(T_i)_{i \in\mathscr{I}}$.
Therefore $f_{t_0}(\gamma,\alpha)$ is microergodic.
\noqed
\end{pf*}

\begin{pf*}{Proof of part~\ref{thmaddmicc}}
First, we consider the case when $(T_i)_{i \in\mathscr{I}}$ has two
different limit points $t_1$ and $t_2$. By part~\ref{thmaddmicb},
$f_{t_1}(\gamma,\alpha)$ and $f_{t_2}(\gamma,\alpha)$ are
microergodic. So, $(\gamma,\alpha)$ is microergodic by the following
lemma (proved in Appendix~\ref{pfsecmiclemmas}):
\noqed\end{pf*}

%
\begin{lemma}\label{lemaddmic02}
Assume there exists $t_1 \ne t_2$ such that both $f_{t_1}(\gamma
_1,\alpha_1) = f_{t_1}(\gamma_2,\alpha_2)$ and $f_{t_2}(\gamma
_1,\alpha_1) = f_{t_2}(\gamma_2,\alpha_2)$. Then $(\gamma_1,\alpha
_1) = (\gamma_2, \alpha_2)$.
\end{lemma}

We now turn to the case when $(T_i)_{i \in\mathscr{I}}$ has only one
limit point $t_0$.
We already know that $f_{t_0}(\gamma,\alpha)$ is microergodic, so we
may assume
that $f_{t_0}(\gamma_1,\alpha_1)=f_{t_0}(\gamma_2,\alpha_2)$, that
is, $g(t_0)=0$.
Denote $\mathscr{I}_{(t_0,\infty)} = \bigcup^\infty_{n=1}{\mathscr
{I}^{\mathbb{T}_n}_{(t_0,\infty)}}$ and $\mathscr{I}_{[0,t_0]} =
\bigcup^\infty_{n=1}{\mathscr{I}^{\mathbb{T}_n}_{[0,t_0]}}$. The
condition in~\ref{thmaddmicc} implies that
$\sum_{i \in\mathscr{I}_{(t_0,\infty)}}(T_i - t_0)^2=\infty$ or
$\sum_{i \in\mathscr{I}_{[0,t_0]}}(T_i - t_0)^2=\infty$
or both. We now use Lemma~\ref{lemaddmic03} (Appendix \ref
{pfsecmicic}) to select, for each $n$, a large set $\mathscr{C}_n$
of independent contrasts such that
$\lim_{n} \sum_{C \in\mathscr{C}_n}{(T_C - t_0)^2} = \infty$.
Again, by (\ref{eqaddmic01}) we have
$r(\mathbb{T}_n) \geq r(\mathscr{C}_n) = \sum_{C \in\mathscr{C}_n} g(T_C)$,
which we approximate below.
If $t_0 = 0$, by Taylor expansion there exists $c(\alpha,T^*)$ such
that $|e^{-2 \alpha x} - 1 + 2 \alpha x - 2 \alpha^2 x^2 + \frac
{4}{3} \alpha^3 x^3| \leq c(\alpha,T^*) x^4$ for every $x$ satisfying
$|x| < T^*$.
Similarly, if $t_0 > 0$ there exists $c(\alpha,T^*)$ such that
$|e^{-2 \alpha x} - 1 + 2 \alpha e^{-2 \alpha t_0} (x-t_0) - 2 \alpha
^2 e^{-2 \alpha t_0} (x-t_0)^2| \leq c(\alpha,T^*) x^3$ for every $x$
satisfying $|x-t_0| < T^*$.
In both cases, we can then write
\[
r(\mathscr{C}_n) = \frac{1}{2}\sum
_{C \in\mathscr{C}_n} {\bigl(h_{t_0}(\alpha_1) -
h_{t_0}(\alpha_2)\bigr)}^2 (T_C-t_0)^2
+ o(T_C-t_0)^2,
\]
where $o(T_C-t_0)^2$ is uniform in $n$, $h_0(\alpha)=\alpha$ and
$h_t(\alpha) = 2\alpha e^{-2 \alpha t}/(1-e^{-2\alpha t})$ for $t>0$.
Therefore $r(\mathscr{C}_n) \to\infty$
unless $(\gamma_1, \alpha_1) = (\gamma_2, \alpha_2)$.
Hence $(\gamma, \alpha)$ is microergodic.
\end{pf}

\begin{figure}

\includegraphics{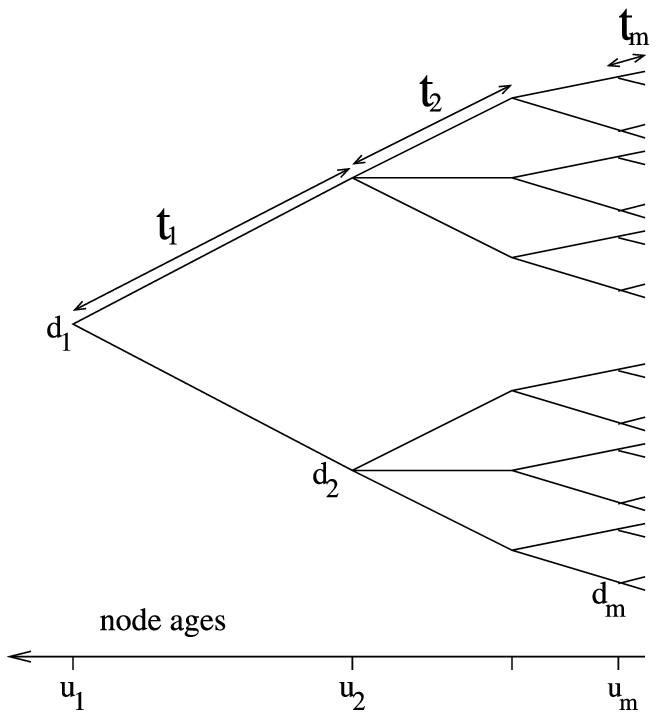}

\caption{Symmetric trees with $m=4$ levels.}
\label{figfig8and4}
\end{figure}

Theorem~\ref{thmaddmic} part~\ref{thmaddmicc} gives a very
general sufficient condition ensuring the microergodicity of $(\gamma,
\alpha)$. Unfortunately, it is not a necessary condition in general.
To prove so, we consider the particular case when $\mathbb{T}_n$ is a
symmetric tree, that is, a tree in which
each internal node is the parent of subtrees of identical shapes
(see Figure~\ref{figfig8and4}). We give below 3 examples in which
%
%
$(T_i)_{i \in\mathscr{I}}$ has only one limit point $t_0$, and the
condition in Theorem~\ref{thmaddmic} part~\ref{thmaddmicc} is violated.
Two examples illustrate the nonmicroergodicity of $(\gamma,\alpha)$,
one in which $t_0 > 0$ and one in which $t_0 = 0$.
In the last example the condition in~\ref{thmaddmicc} is
violated, yet $(\gamma, \alpha)$ is microergodic.

%
\begin{theorem}\label{thmaddmic05}
Consider the OU model (\ref{eqrandomrootV}) on symmetric trees with
$m$ levels
and whose internal nodes at level $i$ have $d_i$ descendants along
branches of length~$t_i$.
\renewcommand\theenumi{(\alph{enumi})}
\renewcommand\labelenumi{\theenumi}
\begin{enumerate}
\item\label{thmaddmic05a}\textit{Increasing node degrees}.
Consider a nested sequence of symmetric trees with a fixed
number of levels $m$ and fixed branch lengths $t_1, t_2,\ldots, t_m$.
Assume that the number of descendants $d_m$ at the last level
goes to infinity, but all other $d_1,\ldots, d_{m-1}$ are fixed,
so that $t_m > 0$ is the only limit point of $(T_i)_{i \in\mathscr{I}}$.
Then $(\gamma,\alpha)$ is not microergodic.

\item\label{thmaddmic05b}\textit{Dense sampling near the tips}, \textit
{or at
distance $t_0$ from
the tips}.
Consider a nested sequence of symmetric trees with a growing number of
levels $m$, $d_k = d$ descendants at all levels $k \geq1$
and such that the age of nodes at level $k$
is $u_k = q^k + t_0$ for some $0 < q < 1$.
Suppose that $dq^2 <1$, to guarantee the violation of the
condition in Theorem~\ref{thmaddmic}\ref{thmaddmicc}:
\begin{enumerate}[(ii)]
\item[(i)] If $t_0=0$, then $(\gamma,\alpha)$ is not microergodic.
\item[(ii)] If $t_0 > 0$, then $(\gamma,\alpha)$ is microergodic.
\end{enumerate}
\end{enumerate}
\end{theorem}

We discuss here the key ingredients of the proof.
The technical details are provided in Appendix~\ref{pfsecmiclemmas}.
Note that node ages are counted with multiplicity.
Here $u_i$ is the age of the $d_1 \lldots d_{i-1}$ internal nodes at
level $i$,
with multiplicity $d_i-1$ for each.
Hence in part~\ref{thmaddmic05a} $u_m=t_m$ is the only limit point.
For a symmetric tree, the eigenvalues of the covariance matrix
are $\gamma\lambda_k(\alpha)$ with multiplicity $d_1 \lldots d_{k-1}
(d_k - 1)$, where
$\lambda_k(\alpha) = \sum_{i=k}^m d_{i+1} \lldots d_m (e^{- 2
\alpha u_{i+1}} - e^{- 2 \alpha u_i})$
(Appendix~\ref{appendix01}).
In~\ref{thmaddmic05a}, only the multiplicity of the smallest eigenvalue
increases to infinity when the tree grows. If
$(\gamma_1, \alpha_1)$ and $(\gamma_2, \alpha_2)$ share the same smallest
eigenvalue, that is, if
$\gamma_1 \lambda_m(\alpha_1) = \gamma_2 \lambda_m(\alpha_2)$,
then insufficient information is gained to distinguish between
$P_{(\gamma_1, \alpha_1)}$ and $P_{(\gamma_2, \alpha_2)}$ when the
tree grows.
In~\ref{thmaddmic05b}, the eigenvalue with the largest
multiplicity is
also the smallest, $\gamma\lambda_m(\alpha) = \gamma(1 - e^{- 2
\alpha u_m})$.
It converges to $0$ when $t_0 = 0$ and to
$\gamma(1 - e^{- 2 \alpha t_0}) > 0$ when $t_0 > 0$, yielding too little
information in (i) when $t_0=0$, but more information to distinguish
between $P_{(\gamma_1, \alpha_1)}$ and $P_{(\gamma_2, \alpha_2)}$
in (ii) when $t_0 > 0$.

\section{Different convergence rates of ML estimators for different
microergodic parameters} \label{secrates}
Section~\ref{secmic} suggests that the different parameters may not be
estimated at the same rate. Indeed, if $t_0$ is the only limit point of
internal node ages, then Theorem~\ref{thmaddmic} shows that
$f_{t_0}(\gamma,\alpha)$ is microergodic regardless of whether condition
in~\ref{thmaddmicc} is satisfied or not. Therefore, the ML or REML estimate
of $f_{t_0}(\gamma,\alpha)$ is expected to converge to the true value
at a faster
rate than the estimate of other parameters. In particular, for $t_0 = 0$
the ML estimate of $\sigma^2$ is expected to converge at a faster rate than
that of~$\alpha$, which might not even be consistent. Here we identify
cases with
unequal convergence rates both theoretically and empirically.

\subsection{\texorpdfstring{Faster convergence of the REML estimator for $f_{t_0}(\gamma,\alpha)$ than for $\alpha$ and~$\gamma$}
{Faster convergence of the REML estimator for f t0(gamma, alpha) than for alpha and gamma}} 

We focus here on the symmetric tree growth model from Theorem~\ref{thmaddmic05} part~\ref{thmaddmic05a} with nodes of
increasing degrees, but we consider here the case when
$\tilde n = n/d_m = d_1\lldots d_{m-1}$ increases indefinitely
to ensure the microergodicity of $\gamma$ and $\alpha$.
We show that the REML estimator of $(\gamma,\alpha)$ is consistent
and asymptotic normally distributed. We further show that
$f_{t_m}(\gamma,\alpha)$,
which is microergodic regardless of the growth of $\tilde n$,
is estimated at a faster rate than $\alpha$ or $\gamma$, which
have stronger requirements to be microergodic.\looseness=-1


%
\begin{theorem}\label{thmreml01}
Consider the asymptotic growth model from above
with OU model (\ref{eqrandomrootV}).
Denote $\nu= \gamma(1 - e^{- 2 \alpha t_m})$. Then the REML estimator
$(\hat{\nu}, \hat{\alpha})$ is consistent and
\[
\pmatrix{\sqrt{n}(\hat{\nu} - \nu)
\cr
\sqrt{\tilde{n}}(\hat{\alpha} - \alpha)}
\stackrel{d} {\rightarrow} N \left( \mathbf{0}, \pmatrix{ 8\nu^2 & 0
\cr
0 & v_{\alpha} } \right).
\]
Moreover, if $n/\tilde n=d_m$ converges to infinity,
then $\sqrt{\tilde n}(\hat\gamma-\gamma,\hat\alpha-\alpha)^{t}$
converges to a centered normal distribution and the asymptotic
correlation between $\log\hat\gamma$ and $\log(1 - e^{- 2 \hat
\alpha t_m})$
is $-1$.\vadjust{\goodbreak}
\end{theorem}

The proof in Appendix~\ref{pfsecrates} gives the expression for
$v_\alpha$.
With increasing node degrees at $m$ levels,
the age of nodes at the last level $t_m$ is the only limit
point of $(T_i)_{i \in\mathscr{I}}$ if $\tilde n$ is bounded.
The growth of $\tilde n$ ensures at least 2 limit points and
the consistency of all parameters.
Our results show that the rate of convergence
is $\tilde n^{-1/2}$ for both $\hat\alpha$ and $\hat\gamma$.
However, only one limit point ($t_m$) is required for the consistent
estimation of $\nu=f_{t_m}(\gamma,\alpha)$, which is microergodic
regardless of $\tilde n$. Accordingly, the convergence rate of
$\hat\nu$ is $n^{-1/2}$, which can be much faster than $\tilde n^{-1/2}$.

\subsection{Simulations on a very large real tree} \label{secsim}

In this section we use simulations to investigate the properties of the
MLE of the OU parameters on a real tree, comprising 4507 mammal species
from \citet{binindaEmonds-etal07}.
%
%
\begin{figure}

\includegraphics{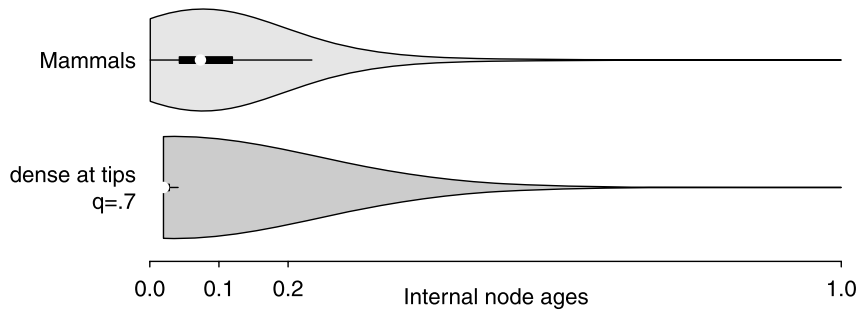}

\caption{Distribution of node ages in the mammal tree (top) and in a
symmetric tree (bottom) of similar size $n=2^{12}=4096$ with $d=2$ at
each level, levels being added near the tips at ages $q^m$.
The value $q=0.7\approx2^{-1/2}$ is the largest at which $\alpha$ and
$\gamma$ are not microergodic.}
\label{fignodedistribution}
\end{figure}
Figure~\ref{fignodedistribution} shows the distribution of node ages for
this tree, and for a symmetric tree with dense sampling near the tips
described in Theorem~\ref{thmaddmic05}\ref{thmaddmic05b},
on which $\alpha$ and $\gamma$ are not microergodic.
Both distributions show a high density of very young nodes.
Under the symmetric tree asymptotics with $0$ as the only limit point,
$\sigma^2$ is microergodic while $(\gamma,\alpha)$ might not be.
Note that this is also the behavior under spatial infill asymptotics in
dimension $d \leq3$.
For real trees like this mammal tree, therefore, we expect the MLE of
$\sigma^2$ to converge quickly, and the MLE of $\alpha$ to converge
more slowly or not at all.
%
%
\begin{figure}

\includegraphics{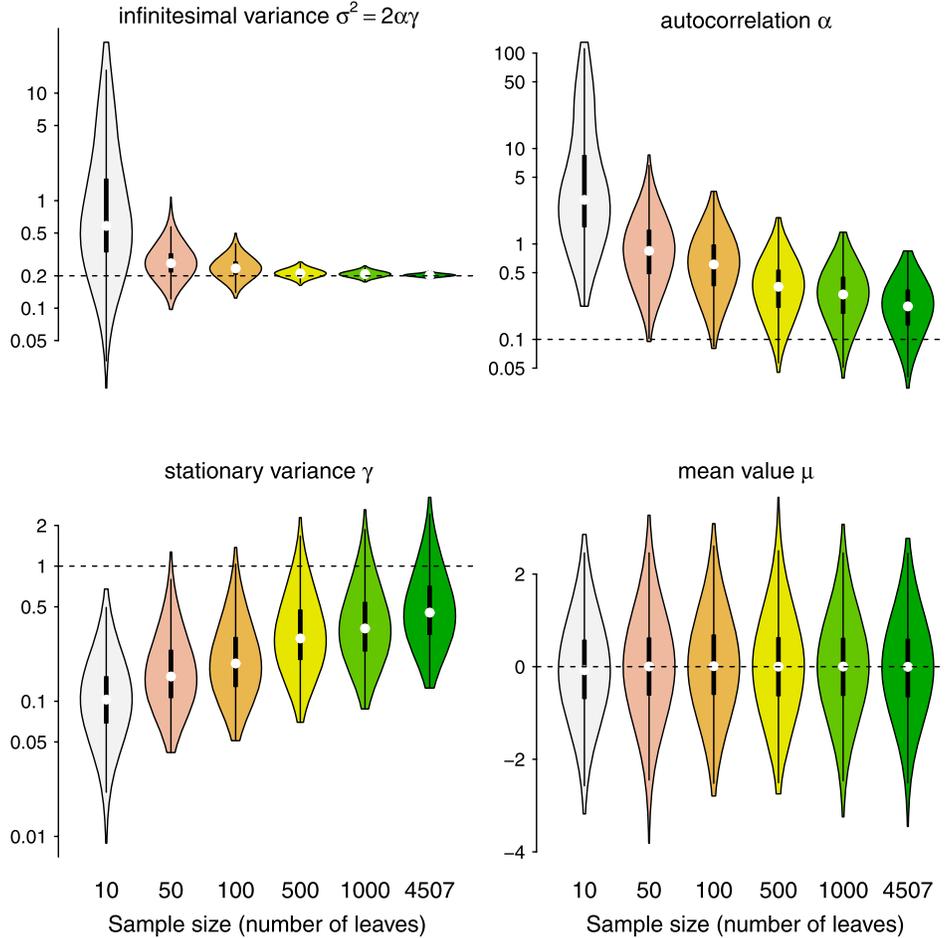}

\caption{Violin plots showing the distribution of the MLE of $\mu$,
$\gamma$, $\alpha$ and $\sigma^2 = 2 \gamma\alpha$ on trees
subsampled from the mammal phylogeny in Bininda-Emonds
et~al. (\citeyear{binindaEmonds-etal07})
with 2000 simulations at each sample size. The true values were $\mu=
0$, $\gamma= 1$, $\alpha= 0.1$ and $\sigma^2 = 0.2$.}
\label{figmammalTreesimulation}
\end{figure}
For various sample sizes from 10 to 4507 (full tree), we simulated data
from the OU model with $\mu= 0$, $\gamma= 1$ and $\alpha= 0.1$, so
$\sigma^2 = 0.2$. We created 20 sequences of six nested trees from
4507 to 10 leaves by randomly selecting subsets of leaves, conditional
on the root being the only common ancestor of the selected leaves to
guarantee that all trees have the same height. Trees were all rescaled
by the same factor to have height 1. For each tree, we simulated 100
data sets and computed the MLEs $\hat{\mu}$, $\hat{\gamma}$ and
$\hat{\alpha}$.
As expected, these simulations show that $\hat\sigma^2$ converges quickly
to the true value while $\hat\alpha$ and $\hat\gamma$ do not
(Figure~\ref{figmammalTreesimulation}).
A strong bias is apparent for $\hat\gamma$ and $\hat\alpha$ even at
the largest sample size (4507). Moreover, the correlation between $\log
\hat{\alpha}$ and $\log\hat{\gamma}$ converges very fast to $-1$
(Table~\ref{tab1}). Also, the lower bound for the variance of $\hat
\mu$ is very close to the true variance (Table~\ref{tab1}).
Therefore, this lower bound can be useful in practice at finite sample sizes.

\section{Discussion}

We considered an Ornstein--Uhlenbeck model of hierarchical
autocorrelation and
showed that the location parameter, here the mean $\mu$, is not
microergodic. We provided the lower bound for the variance of its ML estimator.
In practice, these results could have important implications
when scientists use OU hierarchical autocorrelation
to detect a location shift, that is, a change in $\mu$
along a branch of the tree
[e.g., \citet
{butlerKing04,lavin2008morphometrics,monteiro2011evolutionary}].
Often times, the OU model is used with multiple adaptive optima whose
placements on the tree are not fully known.
Our results suggest that the power to detect such shifts may be low and
mostly influenced by the effect size rather than by the sample size.
An open question is whether the location of such shifts on the tree
can be identified consistently with a growing number of tips.

%
%
\begin{table}
\caption{Correlation between $\log\hat{\alpha}$ and $\log\hat
{\gamma}$
and variance of $\hat{\mu}$ from simulations.
Last line: value of theoretical bound (\protect\ref{eqninc01}) for
$\operatorname{var}(\hat{\mu})$, averaged over 20 simulation subtrees}
\label{tab1}
\begin{tabular*}{\tablewidth}{@{\extracolsep{\fill}}ld{2.4}d{2.4}d{2.4}d{2.4}d{2.4}d{2.4}@{}}
\hline
\textbf{Sample size} & \multicolumn{1}{c}{\textbf{10}} & \multicolumn{1}{c}{\textbf{50}}
& \multicolumn{1}{c}{\textbf{100}} & \multicolumn{1}{c}{\textbf{500}}
& \multicolumn{1}{c}{\textbf{1000}} & \multicolumn{1}{c@{}}{\textbf{4507}}\\
\hline
$\operatorname{cor}(\log\hat\alpha,\log\hat\gamma)$
& -0.44 & -0.927 & -0.9674 & -0.9938 & -0.9971 & -0.9993\\
$\operatorname{var}(\hat{\mu})$ & 0.9007 & 0.8455 & 0.8499 & 0.8853 &
0.8789 & 0.8851\\
Lower bound (\ref{eqninc01}) & 0.8517 & 0.8472 & 0.8469 & 0.8468 &
0.8468 & 0.8468\\
\hline
\end{tabular*}
\end{table}

We provide a general sufficient condition for the covariance parameters
to be microergodic.
Properties of infill asymptotics were recovered when 0 is the only
limit point of internal node ages, that is, when new nodes were added
closer and closer to already existing tips. In this case, $\sigma^2$
is necessarily microergodic. This asymptotics can be appropriate for
coalescent trees or when many species diverged recently from a moderate
number of genera.
We assume here the idealized situation with no error in the tree
structure (topology and branch lengths) and no data measurement error,
leaving this for future work. With measurement error, the covariance
matrix becomes $\gamma\mathbf{V}_\alpha+ \sigma^2_e \mathbf{I}$. The
error variance $\sigma^2_e$ is called a nugget effect in spatial
statistics. Measurement error with tree-structured correlation is
rarely accounted for in applications; but see \citet{ives2007within}.

For a general tree growth model, by using independent contrasts we can
construct a consistent estimator for $f_{t_0}(\gamma, \alpha)$ where
$t_0$ is any limit point of $(T_i)_{i \in\mathscr{I}}$. If $(T_i)_{i
\in\mathscr{I}}$ has at least two limit points, then by Lemma \ref
{lemaddmic02}, we can construct a consistent estimator for $(\gamma,
\alpha)$.
This proposed estimator is based on a restricted set of well-chosen
contrasts, but it uses fewer contrasts and thus less information than
the conventional REML estimator. We conjecture that if $(\gamma,
\alpha)$ is microergodic, the REML estimator of $(\gamma, \alpha)$
is also consistent and asymptotically normal.

The microergodicity results suggest that parameters may not all be
estimated at the same rate.
Indeed, we show that the REML of $\alpha$ converges at a slower rate than
$n^{-1/2}$ under a symmetric tree asymptotic framework.
Similarly, our simulations suggest that the mammalian tree with 4507 species
shares features similar to those under infill asymptotics
(in low dimension) and under dense sampling near the\vadjust{\goodbreak}
tips of symmetric trees, where $\sigma^2$ can be consistently
estimated but
$\alpha$ and $\gamma$ cannot.
On the real tree, the MLE of $\sigma^2$ converges quickly to the true
value while that of $\alpha$ and $\gamma$ do not.
This behavior may explain a lack of power to discriminate between a
model of neutral evolution ($\alpha= 0$) versus a model with natural
selection ($\alpha\ne0$), as observed in \citet{cooperPurvis10}. It
would be interesting to know if most real trees share the ``dense tip''
asymptotic behavior,
or how frequently a ``dense root'' asymptotic is applicable instead.
Our results point to the distribution on node ages
as indicative of the most appropriate asymptotic regime.

%
\begin{appendix}
\section{Spectral decomposition of the OU covariance matrix on
symmetric trees} \label{appendix01}
We consider here symmetric trees (Figure~\ref{figfig8and4})
with $m$ levels of internal nodes, the root being at level $1$.
Each node at level $k$ is connected to $d_k\geq2$ children
by branches of length $t_k$. The age of nodes at level $k$
is then $u_k = t_k+\cdots+t_m$.
Under the OU model (\ref{eqrandomrootV}),
the correlation matrix $\mathbf{V}_\alpha$
is identical to that obtained under a BM
model along a tree with an extra branch extending from the root
and with transformed branch lengths $\mathbf{t}^{\mathrm{BM}}$,
\[
t^{\mathrm{BM}}_k(\alpha) = \cases{ 1 - e^{- 2 \alpha t_m}, &\quad for
$k=m$,
\cr
e^{-2\alpha u_{k+1}} - e^{-2\alpha u_k}, &\quad$1\leq k \leq m-1$,
\cr
e^{-2 \alpha u_1}, &\quad$k=0$ (extra root branch).}
\]
Therefore, we can derive the eigen-decomposition of
$\mathbf{V}_{\alpha}(\mathbf{t}) = \mathbf{V}^{\mathrm{BM}}_{\alpha
}(\mathbf{t}^{\mathrm{BM}})$
as done in \citet{ane08}.
The eigenvalues, from greatest to smallest, are
\[
\lambda_k = n \sum_{i=k}^m{
\frac{t^{\mathrm{BM}}_i(\alpha)}{d_1
\lldots d_i}} = \sum_{i=k}^m
d_{i+1} \lldots d_m \bigl(e^{- 2 \alpha u_{i+1}} -
e^{- 2
\alpha u_i}\bigr) 
\]
with multiplicity $d_1 \lldots d_{k-1} (d_k - 1)$, for $k = 0,\ldots,m$
and $u_{m+1}$ set to $0$ and $u_0$ to~$\infty$.
Furthermore, \citet{ane08} showed that
the eigenvectors of $\mathbf{V}^{\mathrm{BM}}_{\alpha}$ are independent of the tree's
branch lengths, which implies here that the eigenvectors of
$\mathbf{V}_{\alpha}$ are independent of $\alpha$.
Each eigenvector corresponding to $\lambda_k(\alpha)$ represents a contrast
between the descendants of a node at level $k$. One exception is the eigenvector
associated with the extra root branch and largest eigenvalue $\lambda_0$.
This eigenvector is $\one$ and has multiplicity $1$.

\section{Supporting lemmas and technical proofs} \label{appendix02}

\subsection{Procedures for choosing independent contrasts}
\label{pfsecmicic}
%
%
\begin{lemma}\label{lemaddmic05}
Let $\mathbb{T}$ be an ultrametric tree. For every $a < b$, we can
choose a set of independent contrasts $\mathscr{C}$ with respect to
some of the internal nodes in $\mathscr{I}^{\mathbb{T}}_{(a,b)}$ such that
$|\mathscr{C}| \geq\frac{1}{2} |\mathscr{I}^{\mathbb{T}}_{(a,b)}|$.
\end{lemma}

\begin{pf}
We choose contrasts as follows, starting with
$\mathscr{C} = \varnothing$ and $\mathbb{T}_0 = \mathbb{T}$.
At step $n$, we choose an\vadjust{\goodbreak} internal node $i_n \in\mathscr{I}^{\mathbb
{T}_{n-1}}_{(a,b)}$
of minimum age, and a path $p_{i_n}$ connecting any two tips having
$i_n$ as their common ancestor. We update $\mathscr{C} = \mathscr{C}
\cup\{ C^{p_{i_n}}_{i_n} \}$ and obtain tree $\mathbb{T}_n$ from
$\mathbb{T}_{n-1}$ by dropping all descendants of $i_n$. We stop when
$\mathscr{I}^{\mathbb{T}_{n}}_{(a,b)} = \varnothing$.
The procedure guarantees that the paths do not intersect, hence the
contrasts are independent. Furthermore, $\mathscr{I}^{\mathbb
{T}_n}_{(a,b)} = \mathscr{I}^{\mathbb{T}_{n-1}}_{(a,b)} \setminus\{
i_n, i'_n \}$ where $i'_n$ is the parent of $i_n$, so $|\mathscr{C}|
\geq|\mathscr{I}^{\mathbb{T}}_{(a,b)}|/2$.
\end{pf}

%
\begin{lemma}\label{lemaddmic03}
Let $\mathbb{T}$ be an ultrametric tree of height $T$. For all $t \in[0,T]$:
\begin{longlist}[(a)]
\item[(a)] There exists a set of independent contrasts $\mathscr{C}$ with
respect to nodes in $\mathscr{I}^{\mathbb{T}}_{[0,t]}$ such that
$\sum_{C \in\mathscr{C}}{(T_C - t)^2} \geq\frac{1}{2} \sum_{i \in
\mathscr{I}^{\mathbb{T}}_{[0,t]}}{(T_i - t)^2}$.
\item[(b)] There exists a set of independent contrasts $\mathscr{C}$ with
respect to nodes in $\mathscr{I}^{\mathbb{T}}_{(t,\infty)}$ such that
$
\sum_{C \in\mathscr{C}}{(T_C - t)^2} \geq\frac{1}{4} [(T -
t)^2 + \sum_{i \in\mathscr{I}^{\mathbb{T}}_{(t,\infty)}}{(T_i -
t)^2} ].
$
\end{longlist}
\end{lemma}

\textit{Proof of Lemma~\ref{lemaddmic03}.\quad}
(a)
The procedure in the proof of Lemma~\ref{lemaddmic05} gives us a
desired set of contrasts. Indeed, let $(i_k)_{k=1}^m$ be the chosen set
of nodes and $(i'_k)_{k=1}^m$ be their parents. Then $\mathscr
{I}^{\mathbb{T}}_{[0,t]} \subset\bigcup_{k=1}^m \{ i_k, i'_k \}$, hence
\[
\sum_{i \in\mathscr{I}^{\mathbb{T}}_{[0,t]}} (T_i - t)^2 \leq
\sum_{k=1}^m (T_{i_k}-t)^2
+ (T_{i'_k}-t)^2 \leq2 \sum_{k=1}^m
(T_{i_k}-t)^2 = 2 \sum_{C \in\mathscr{C}}
(T_C-t)^2.
\]

%
%
\begin{figure}[b]

\includegraphics{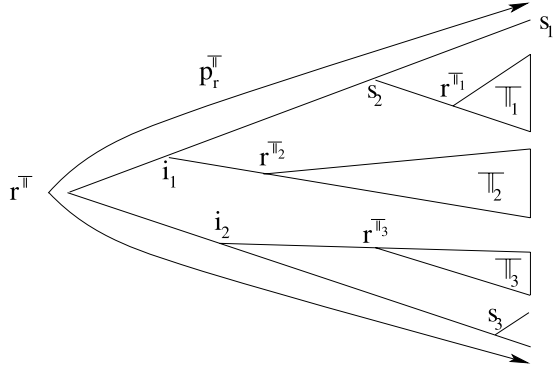}

\caption{Recursive construction of independent contrasts, taken with
respect to the root at each step.}
\label{figfig7}
\end{figure}

(b)
Contrasts are chosen by induction, starting with $\mathscr{C} =
\varnothing$. Let $r^{\mathbb{T}}$ be the root of $\mathbb{T}$. If
$r^{\mathbb{T}} \notin\mathscr{I}^{\mathbb{T}}_{(t,T]}$, then we
stop; else we update $\mathscr{C} = \mathscr{C} \cup\{
C^{p_r^{\mathbb{T}}}_{r^{\mathbb{T}}} \}$ where the path
$p_r^{\mathbb{T}}$ is chosen carefully as follows. From each child of
the root, the path descends toward the tips. Each time an internal node
is encountered, a decision needs to be made to either go left or right.
Of the two children of the internal node, the path is connected to the
youngest (Figure~\ref{figfig7}).
We then remove from $\mathbb{T}$ the path $p_r^{\mathbb{T}}$ and the
edges connected to it. What is left is a forest, a set of subtrees of
$\mathbb{T}$, one which we repeat the procedure, recursively
extracting one path and its corresponding contrast from each subtree.

We now prove by induction that this procedure gives us a desired set of
contrasts. This is easy to see for $\leq3$ tips. Assume that it is
true for every tree with \mbox{$\leq$}$m$ tips, and that $\mathbb{T}$ has
$m+1$ tips. Let $i_1$ and $i_2$ be the two children of $r^{\mathbb
{T}}$. Let $(\mathbb{T}_k)_{k=1}^l$ be the subtrees obtained after
removing $p_r^{\mathbb{T}}$ and the edges connected to it, and such
that $r^{\mathbb{T}_k} \in\mathscr{I}^{\mathbb{T}}_{(t,T]}$. Let
$s_k$ be the sibling of $r^{\mathbb{T}_k}$ in $\mathbb{T}$
($s_k$ could be a leaf). By construction, $T_{s_k} \leq T_{r^{\mathbb
{T}_k}}$. Let $\mathscr{C}_k$ be the set of contrasts obtained from
$\mathbb{T}_k$. We have $\mathscr{I}^{\mathbb{T}}_{(t,T]} \subset\{
r^{\mathbb{T}}, i_1, i_2\} \bigcup_{k=1}^l \mathscr{I}^{\mathbb
{T}_k}_{(t,T]} \cup\{ s_k \}$ and $\mathscr{C} = \{ r^{\mathbb
{T}} \} \bigcup_{k=1}^l \mathscr{C}_k$. Therefore,
\begin{eqnarray*}
4 \sum_{C \in\mathscr{C}}{(T_C-t)^2} &=&
4(T_{r^{\mathbb{T}}}-t)^2 + 4 \sum_{k=1}^l
{\sum_{C \in\mathscr{C}_k}{(T_C-t)^2}}
\\
&\geq& 2(T_{r^{\mathbb{T}}}-t)^2 + \bigl(\max\{ T_{i_1},t \}-t
\bigr)^2 + \bigl(\max\{ T_{i_2},t \}-t\bigr)^2
\\
&&{} + \sum_{k=1}^l \biggl\{
(T_{r^{\mathbb{T}_k}} - t)^2 + \sum_{i \in
\mathscr{I}^{\mathbb{T}_k}_{(t,T]}}
{(T_i-t)^2} \biggr\} \\
&\geq&(T_{r^{\mathbb{T}}}-t)^2
+ \sum_{i \in\mathscr{I}^{\mathbb
{T}}_{(t,T]}}{(T_i - t)^2}.
\end{eqnarray*}

\subsection{\texorpdfstring{Technical proofs for Section \protect\ref{secmic}}{Technical proofs for Section 2}}
\label{pfsecmiclemmas}

\subsubsection*{Counter example for Theorem \protect\ref{thminc01}
on nonultrametric trees}

Let $a = e^{-\alpha t_1}$ and $b = e^{-\alpha t_2}$.
It is easy to see that $\mathbf{V}_\alpha$ can be expressed in
terms of the $n/2\times n/2$ identity matrix $\mathbf{I}$ as
$
\mathbf{V}_\alpha= \operatorname{diag}((1-a^2)\mathbf{I}, (1-b^2)\mathbf{I})
+ (a\one^{t}, b\one^{t})^{t}
(a\one^{t}, b\one^{t})$.
We then get $\mathbf{V}_\alpha^{-1}$ using Woodbury's formula, then
$
\one^{t} \mathbf{V}^{-1}_\alpha\one= n ( \frac{1}{1-a^2} +
\frac{1}{1-b^2} + \frac{n(a-b)^2}{(1-a^2)(1-b^2)} )
/ ( 1 + \frac{n a^2}{1-a^2} + \frac{n b^2}{1-b^2} ).
$
If $t_1\neq t_2$, then $a \ne b$ and
$\operatorname{var}(\hat{\mu}) = (\one^{t} V^{-1}_\alpha\one)^{-1}$
goes to 0 as claimed.\vspace*{9pt}

\textit{Proof of Lemma~\ref{leminc02}.\quad}
We will first prove
$\mathbf{V}_{\alpha} \geq e^{-2 \alpha T} \mathbf{J}_n$ by induction
on the
number of tips, where $\mathbf{J}_n = \one\one^{t}$.
Clearly, this is true for trees with a single tip.
Now consider a tree with $n$ tips, and consider its $k$ subtrees
obtained by removing the $k$ branches stemming from the root.
Let $T_1,\ldots, T_k$ be the heights of these subtrees, that is,
the age of their roots.
Their number of tips $n_1,\ldots, n_{k}$ is at most $n-1$.
So by induction, the covariance matrices
$\mathbf{V}^{(1)}_\alpha,\ldots, \mathbf{V}^{(k)}_\alpha$
associated with these subtrees
must satisfy
$\mathbf{V}^{(i)}_\alpha\geq e^{-2 \alpha T_i} \mathbf{J}_{n_i}$. Therefore
$\mathbf{V}_{\alpha} - e^{-2 \alpha T} \mathbf{J}_n \geq
\operatorname{diag}(\mathbf{V}^{(i)}_\alpha- e^{-2 \alpha T_i} \mathbf
{J}_{n_i}) \geq0$
is true for all trees.
Now we use the definition\vadjust{\goodbreak} of $t$ and go a step further using that
$
\mathbf{V}^{(i)}_\alpha- e^{-2 \alpha T} \mathbf{J}_{n_i} \geq
(e^{-2 \alpha T_i} - e^{-2 \alpha T}) \mathbf{J}_{n_i}
\geq(e^{-2 \alpha(T-t)} - e^{-2 \alpha T}) \mathbf{J}_{n_i}
$
for all $i=1,\ldots,k$. This implies that
$
\mathbf{V}_{\alpha}-e^{-2\alpha T}\mathbf{J}_n
\geq(e^{-2\alpha(T-t)} - e^{-2\alpha T})
\operatorname{diag}( \mathbf{J}_{n_1},\ldots, \mathbf{J}_{n_{k}})
\geq\frac{1}{k}(e^{-2 \alpha(T-t)} - e^{-2 \alpha T}) \mathbf{J}_n,
$
from which Lemma~\ref{leminc02} follows easily.\vspace*{9pt}

\textit{Proof of upper bound \textup{(\ref{lemaddmic04})}.\quad}
Assume here that $\alpha_1=\alpha_2$ and $\gamma_1=\gamma_2=\gamma$.
Since $(Y_i)_{i=1}^n$ have the same covariance matrix $\gamma{\mathbf V}$
under both distributions $P_{\theta_1}$ and $P_{\theta_2}$, it is
easy to see that
$r(\mathbb{T})=(\mu_1-\mu_2)^2 \one^{t} {\mathbf V}^{-1}\one
/\gamma$
[\citet{hersheyOlsen07}].
The bound $r(\mathbb{T}) \leq(\mu_1 - \mu_2)^2/(\gamma e^{- 2
\alpha T})$,
where $T$ is the age of the root, then follows from Lemma~\ref{leminc02}.\vspace*{9pt}

\textit{Proof of Lemma~\ref{lemaddmic06}.\quad}
First, $\operatorname{var}(C^{p_i}_i) = \operatorname{var}(Y_1^i) +
\operatorname{var}(Y_2^i) - 2 \operatorname{cov}(Y_1^i,Y_2^i) = 2 \gamma
- 2 \gamma e^{-2
\alpha T_i}$.
Second, consider two paths $p_{i_1}$ and $p_{i_2}$ that do not intersect.
Then, the most recent common ancestor of $Y_j^{i_1}$ and $Y_k^{i_2}$
($j,k \in\{ 1,2 \}$) is the most recent common ancestor of internal
nodes $i_1$ and $i_2$. Therefore, the distance from $Y_1^{i_1}$ to
$Y_1^{i_2}$ equals the distance from $Y_2^{i_1}$ to $Y_1^{i_2}$. Hence
$\operatorname{cov}(Y_1^{i_1}, Y_1^{i_2}) = \operatorname{cov}(Y_2^{i_1},
Y_1^{i_2})$. Similarly, $\operatorname{cov}(Y_1^{i_1}, Y_2^{i_2}) =
\operatorname{cov}(Y_2^{i_1}, Y_2^{i_2})$. Therefore
$\operatorname{cov}(C^{p_i}_{i_1},C^{p_i}_{i_2}) = \operatorname
{cov}(Y_1^{i_1} -
Y_2^{i_1}, Y_1^{i_2} - Y_2^{i_2}) = 0$.\vspace*{9pt}

\textit{Proof of Lemma~\ref{lemaddmic02}.\quad}
Define $h_1(x) = (1-e^{-2x \alpha_2})/(1-e^{-2x \alpha_1})$, and
assume $t_1 \ne0$ and $t_2 \ne0$.
From the system of equations, we have
$\gamma_1/\gamma_2 = h_1(t_1)=h_1(t_2)$.
Now $(\log h_1)'(x)/x = h_2(x\alpha_2)-h_2(x\alpha_1)$ where
$h_2(x) =\break x e^{-x}/(1-e^{-x})$ is monotone on $(0, \infty)$.
So $\alpha_1 = \alpha_2$, and $\gamma_1 = \gamma_2$.
If $t_2 = 0$, we make a similar argument because
$h_3(x) = x/(1-e^{-2xt_2}) $ is monotone on $(0,\infty)$.\vspace*{9pt}

\textit{Proof of Theorem \protect\ref{thmaddmic05} part
\protect\ref{thmaddmic05a}.\quad}
Under the symmetric tree growth model,
\[
r(\mathbb{T}) = \frac{1}{2}\sum_{k=1}^{m}{d_1
\lldots d_{k-1}(d_k - 1) \biggl( \frac{\gamma_2\lambda_k (\alpha
_2)}{\gamma_1\lambda_k(\alpha_1)} +
\frac{\gamma_1\lambda_k (\alpha_1)}{\gamma_2\lambda_k(\alpha_2)} - 2
\biggr)} + \biggl( m_{1,n}^2 +
\frac{m_{1,n}^2}{\sigma_{1,n}^2} \biggr).
\]
To show this, we consider
$\mathbf{h}=(Y_{j,n})_{j\leq n}=\gamma_1^{-1/2}\bolds{\Lambda
}^{-1/2}(\alpha_1)\mathbf{P}^{-1}(\mathbf{Y}-\mu_1\one)$,
where $\bolds{\Lambda}(\alpha) = \operatorname{diag}(\lambda_k(\alpha))$
contains the eigenvalues $\lambda_k$ with their multiplicities,
and $\mathbf{P}$ contains the eigenvectors of $\mathbf{V}_\alpha$, which
do not depend of $\alpha$ (Appendix~\ref{appendix01}). Then $\mathbf{h}$
is orthonormal under $P_{\theta_1}$, and orthogonal under $P_{\theta
_2}$ with
variances $(\gamma_2/\gamma_1) \lambda_k(\alpha_2)/\lambda
_k(\alpha_1)$ with
multiplicities $d_1 \lldots d_{k-1}(d_k - 1)$.
Furthermore,
$E_2 \mathbf{h} = (\mu_2 - \mu_1)\gamma_1^{-1/2} \bolds{\Lambda
}^{-1/2}(\alpha_1) \mathbf{P}^{-1} \one$
so that $m_{j,n} = 0$ if $j \geq2$ and
$m_{1,n} = (\mu_2 - \mu_1)/\sqrt{n\gamma_1 \lambda_0(\alpha_1)}$,
from which $r(\mathbb{T})$ follows.

With increasing node degrees at $m$ levels, it is easy to
see that the ratio $\lambda_k(\alpha_1)/\lambda_k(\alpha_2)$
converges to a positive limit for all $k\leq m$.
Under the assumption that $d_k$ is fixed for $k<m$,
the multiplicity of $\lambda_k(\alpha)$ is constant as $n$ grows,
except for $k=m$. 
$r(\mathbb{T}_m)$ is then expressed as
a finite sum where all terms are convergent except for the last term
($k=m$) associated with the smallest eigenvalue
$\lambda_m = 1 - e^{- 2 \alpha t_m}$. This term
is bounded if and only if $\gamma_1(1-e^{-2\alpha_1 t_m})=\gamma
_2(1-e^{-2\alpha_2 t_m})$,\vadjust{\goodbreak}
in which case $r(\mathbb{T}_n)$ converges to a finite value.
Otherwise, $r(\mathbb{T}_n)$ goes to infinity. Hence
$P_{\theta_1}$ and $P_{\theta_2}$ are equivalent if and only if
$\gamma_1(1-e^{-2\alpha_1 t_m})=\gamma_2(1-e^{-2\alpha_2 t_m})$,
which completes the proof.\vspace*{9pt}

\textit{Proof of Theorem~\ref{thmaddmic05} part
\textup{\ref{thmaddmic05b}}.\quad}
We denote here $\lambda_k=\lambda_{k,m}$ to emphasize the dependence
of $m$.
We first consider case (i) when $t_0=0$.
When $d_k = d$ and \mbox{$u_k=q^k$}, the eigenvalues simplify to
%
%
\begin{equation}
\label{eqlambdakm} \frac{\lambda_{k,m}(\alpha)}{d^{m-k}} = \sum
_{j=0}^{m-k-1} \frac{e^{-2 \alpha q^{k+1+j}} - e^{-2 \alpha
q^{k+j}}}{d^j} + \frac{1 - e^{-2 \alpha q^m}}{d^{m-k}}.
\end{equation}
It is then easy to see that for all $\alpha\geq0$ and $k$,
$\lambda_{k,m}(\alpha)/d^{m-k}$
converges to some finite function of $\alpha$ and $k$.
To prove the convergence of $r(\mathbb{T}_m)$ we will need the
following lemma,
which is proved later.
%
%
\begin{lemma}\label{lemmic02}
Let $\gamma_1 \alpha_1 = \gamma_2 \alpha_2$, that is, $\sigma_1^2
=\sigma_2^2$.
Then there exists $K$, $c$ and $C$ which depend only on $\alpha_1,
\alpha_2, d$ and $q$
such that for all $m>k \geq K$,
\[
c q^{2k} \leq\frac{\gamma_2\lambda_{k,m} (\alpha_2)}{\gamma_1\lambda
_{k,m}(\alpha_1)} + \frac{\gamma_1\lambda_{k,m} (\alpha_1)}{\gamma
_2\lambda
_{k,m}(\alpha_2)} - 2 \leq C
q^{2k}.
\]
\end{lemma}
Because $\gamma\alpha$ is microergodic [Theorem~\ref{thmaddmic}
part~\ref{thmaddmicb}], we can assume $\gamma_1 \alpha_1 = \gamma
_2 \alpha_2$.
Lemma~\ref{lemmic02} implies that the first sum
in the expression of $r(\mathbb{T}_m)$ [from the proof of part (a)]
is bounded above and below by $\sum_{k=K}^m{(d q^2)^k}$ up to some
multiplicative constant,
and so converges to a finite limit because $d q^2 < 1$.
The last term with $m_{1,n}^2$ is always bounded as shown in the
proof of (\ref{lemaddmic04}). This completes the proof.

We now turn to case (ii) with $t_0>0$. To prove that $(\gamma,\alpha
)$ is microergodic, we will show that $P_{\theta_1} \,\bot\, P_{\theta
_2}$ under the restriction $\gamma_1 (1-e^{-2t \alpha_1}) = \gamma_2
(1-e^{-2t \alpha_2})$. To do so, we only need to check the sufficient
condition in (\ref{condRao}). Note that there exits $w > 0$ such that
$d^w q \geq1$. Denote $k_m = [ m/(w+1) ]$ where $[x]$ is a largest
integer smaller than $x$.
The condition in (\ref{condRao}), denoted by $z_m$, can be written as
\begin{eqnarray*}
z_m &=& \sum_{k=1}^{m}
d^{k-1}(d-1) \biggl(\frac{\gamma_1\lambda_{k,m}
(\alpha_1)}{\gamma_2\lambda_{k,m}(\alpha_2)} - 1 \biggr)^2 \geq
d^{k_m-1} \biggl( \frac{\gamma_1\lambda_{k_m,m} (\alpha_1)}{\gamma
_2\lambda_{k_m,m}(\alpha_2)} - 1 \biggr)^2
\\
&\geq& d^{k_m-1} \biggl( \frac{(h_{t_0}(\alpha_1) - h_{t_0}(\alpha
_2))f_{m,1} + O_{\alpha_1,\alpha_2,t_0}(1)q^{k_m}}{1+h_{t_0}(\alpha
_2)f_{m,1} + O_{\alpha_2,t_0}(1)q^{k_m}} \biggr)^2,
\end{eqnarray*}
where $ h_t(\alpha) = \frac{2\alpha e^{-2 \alpha
t}}{1-e^{-2\alpha t}}$ and $ f_{m,1}(q) = \sum_{j=0}^{m-k-1}{ ( \frac
{q}{d} )^j} + \frac{1}{1-q}
( \frac{q}{d} )^{m-k}$. If $(\gamma_1,\alpha_1) \ne(\gamma
_2,\alpha_2)$, then $z_m \to\infty$ because $h_t(\alpha)$ is
monotone in $\alpha$.\vspace*{9pt}

\textit{Proof of Lemma~\ref{lemmic02}.\quad}
We first note that for every $a>0$ there exists $x_a>0$ such that
$e^{-ax} - (1 - ax + a^2 x^2/2) = O(a^3x^3)$
uniformly for all $x$ in $[0,x_a]$.
Therefore there exists $K=K(\alpha,q)$ such that
$
e^{-2 \alpha q ^{k+j+1}} - e^{- 2 \alpha q ^{k+j}} -
2 \alpha q^{k+j}(1-q) + 2 \alpha^2 q^{2k+2j}(1-q^2) 
= q^{3k+3j} O_\alpha(1)
$\vadjust{\goodbreak}
where the $O_\alpha(1)$ term is bounded uniformly in $k+j\geq K$.
We can now combine this with
(\ref{eqlambdakm}),
$ \lambda_{k,m}(\alpha)/d^{m-k} - 2 \alpha(1-q) q^k f_1
+ 2 \alpha^2 (1-q^2)q^{2k} f_2
= q^{3k} f_3 O_\alpha(1)
$
where $f_1, f_2$ and $f_3$ only depend on $q,d,m-k$ and are defined by
$f_1=f_{m,1}(q)$, $f_2=f_{m,1}(q^2)$ and
$f_3=f_{m,1}(q^3)$.
Because the $f$ values are bounded as $m-k$ grows, we get
$
\frac{\lambda_{k,m}(\alpha_1)}{\lambda_{k.m}(\alpha_2)}
= \frac{\alpha_1}{\alpha_2}
(1+(\alpha_2-\alpha_1)(1+q)q^k f_2/f_1 ) + q^{2k} O(1)
$
where the\vspace*{2pt} $O(1)$ term is bounded uniformly in $m>k \geq K$, and
the same formula holds when $\alpha_2$ and $\alpha_1$ are switched.
Lemma~\ref{lemmic02}
then follows immediately because we assume that $\gamma_1\alpha
_1=\gamma_2\alpha_2$.

\subsection{\texorpdfstring{Technical proofs for Section \protect\ref{secrates}}{Technical proofs for Section 3}}
\label{pfsecrates}

\subsubsection*{Criterion for the consistency and asymptotic normality
of REML estimators}
In Appendix~\ref{appendix01}, we showed that $\one$ is an eigenvector of
$\mathbf{V}_\alpha$ for symmetric trees, independently of $\alpha$.
Therefore,\vspace*{1pt} the REML estimator of $(\gamma,\alpha)$ based on
$\mathbf{Y}$ is the ML estimator of $(\gamma,\alpha)$ based on the
transformed data $\tilde{\mathbf{Y}} = \tilde{\mathbf
{P}}^{t}\mathbf{Y}$ where
$\tilde{\mathbf{P}}$ is the matrix of all eigenvectors but $\one$.
$\tilde{\mathbf{Y}}$ is Gaussian centered with variance
$\bolds{\Sigma}_n=\gamma\tilde{\bolds{\Lambda}}$ where
$\tilde{\bolds{\Lambda}}$ is the diagonal
matrix of all eigenvalues of $\mathbf{V}_\alpha$ but $\lambda
_0(\alpha)$.
Following \citet{mardiaMarshall84} and like \citet{cressieLahiri93},
we use a general result from \citet{sweeting80}.
The following conditions, C1--C2,
ensure the consistency and asymptotic normality of the ML estimator
[reworded from \citet{mardiaMarshall84}].
Assume there exists nonrandom continuous symmetric matrices
$\mathbf{A}_n(\bolds{\theta})$ such that:
\begin{enumerate}[(C2)]
\item[(C1)] \mbox{}\hphantom{i}(i) As $n$ goes to infinity $\mathbf{A}_n^{-1}$ converges
to $0$.


\noindent (ii) $\mathbf{A}_n^{-1} {\cal J}_n \mathbf{A}_n^{-1}$ converges in probability
to a positive definite matrix $\mathbf{W}(\bolds{\theta})$,
where ${\cal J}_n$ is the second-order derivative of the negative log likelihood
function $L$.
\item[(C2)] $\bolds{\Sigma}_n$ is twice continuously differentiable
on $\Theta$
with continuous second derivatives.
\end{enumerate}
Under these conditions, the MLE $\hat{\bolds\theta}$ satisfies
$\mathbf{A}_n(\bolds{\theta})(\hat{\bolds\theta}-\bolds
{\theta}) \stackrel{d}{\rightarrow}
N (0,\mathbf{W}(\bolds{\theta})^{-1} )$.
A~standard choice for $\mathbf{A}_n$ is the inverse of the
square-root of the Fisher information matrix $\mathbf{B}_n=\Ef({\cal J}_n)$.
Because (C1)(ii) is usually difficult to verify,
\citet{mardiaMarshall84} suggest using a stronger
$L^2$-convergence condition.
This approach was later taken
by Cressie and Lahiri (\citeyear{cressieLahiri93,cressieLahiri96}).
Unfortunately, their conditions for establishing (C1)
do not hold here, because the largest eigenvalues
and the ratio of the largest to the smallest eigenvalues are both
of order $n$.
In what follows, we will check (C1) for the particular
choice of $\mathbf{A}_n=\mathbf{B}_n^{1/2}$ and
$\mathbf{W}(\bolds{\theta}) = \mathbf{I}$
and where we replace (C1)(ii) by the stronger condition
\begin{enumerate}[(C1)]
\item[(C1)](ii$'$)
$\sum_{i,j,k,l=1,2} b^{k i}b^{l j} \operatorname{tr}(\bolds{\Sigma
}_n{(\bolds{\Sigma}_n^{-1})}_{k j}\bolds{\Sigma}_n{(\bolds
{\Sigma}_n^{-1})}_{l i})$
converges to $0$,
where $b^{i j}$ is the $(i,j)$-element of $\mathbf{B}_n^{-1}$,
and ${(\bolds{\Sigma}_n^{-1})}_{i j}$ is the $(i,j)$-second order
derivative of $\bolds{\Sigma}_n^{-1}$.\vadjust{\goodbreak}
\end{enumerate}
%

\textit{Proof of Theorem~\ref{thmreml01}.\quad}
It is convenient here to re-parametrize the
model using $(\nu,\alpha)$.
The diagonal elements in $\bolds{\Sigma}_n$ are
$\nu\lambda_k(\alpha)/\lambda_m(\alpha)$ with
multiplicity $d_1\lldots d_{k-1}(d_k-1)$.
The smallest is $\nu$ (for $k=m$) with multiplicity $n-\tilde n$,
which is conveniently independent of $\alpha$.
With this parametrization, the inverse of the Fisher
information matrix is the symmetric matrix
\[
\mathbf{B}_n^{-1} = \frac{2}{\det\mathbf{B}_n} \pmatrix{\displaystyle\sum
_{k=1}^{m-1} d_1 \lldots
d_{k-1}(d_k - 1) (\Lambda_{k,m} -
\Lambda_{m,m})^2 & *
\vspace*{2pt}\cr
\displaystyle-{\nu}^{-1}
\sum_{k=1}^{m-1} d_1 \lldots
d_{k-1}(d_k - 1) (\Lambda_{k,m} -
\Lambda_{m,m}) & (n-1)/\nu^2 },
\]
where
$\Lambda_{k,m}=\lambda_{k,m}'/\lambda_{k,m}$,
$\det\mathbf{B}_n = (n-1)^2/(4\nu^2) \operatorname{var}_q(\Lambda
_{K,m}-\Lambda_{m,m})$
and the variance is taken with respect to
$\Pf\{K=k \} = q_{k,n}=d_1 \lldots d_{k-1} (d_k - 1)/(n - 1)$.
When the degree at the last level near the tips $d_m$ becomes large
then $q_{m,n}\sim1$, that is, the distribution $q$ is concentrated
around the high end $K=m$.
It is then useful to express
\[
\det\mathbf{B}_n = \frac{(n-{\tilde n})(\tilde n-1)}{4\nu^2} \Ef_p(
\Lambda_{K,m}-\Lambda_{m,m})^2 +\frac{(\tilde n-1)^2}{4\nu^2}
\operatorname{var}_p(\Lambda_{K,m}-\Lambda_{m,m}),
\]
where the expectation and variance are now taken with respect to
$\Pf\{K=k \}=p_{k,n}=d_1 \lldots d_{k-1} (d_k - 1)/(\tilde
n - 1)$ for $k<m$,
that is, $p_{k,n}=q_{k,n}(n-1)/(\tilde n -1)$.
To verify conditions (C1)(i) and (ii$'$), we will use the following lemmas.

%
\begin{lemma}\label{lemreml01}
$\Lambda_{1,m} < \Lambda_{2,m} < \cdots< \Lambda_{m,m}$.
Moreover for any fixed $T$ and $\alpha>0$,
$\Lambda_{k,m}$ and $\lambda_{k,m}''/\lambda_{k,m}$ are uniformly bounded.
Specifically,
$|\Lambda_{k,m}| \leq\max\{ 2T,1/\alpha\}$ and
$|\lambda_{k,m}''/\lambda_{k,m}| \leq4 \max\{ T^2,T/\alpha\}$.
\end{lemma}

\textit{Proof of Lemma~\ref{lemreml01}.\quad}
Denote $g(\alpha)=b-(b+c)e^{-c\alpha}+ce^{-(b+c)\alpha}$.
It is easy to see that
$g'>0$ 
then $g>0$ for all
$\alpha,b,c > 0$. It follows that
\[
\frac{(a+b)e^{- \alpha b} - a}{1 - e^{- \alpha b}} - \frac
{(a+b+c)e^{- \alpha c} - (a+b)}{1 - e ^{- \alpha c}} > 0\qquad\forall
a \in
\mathbb{R},
\alpha,b,c > 0.
\]
Now let $a_i>0$ for $i = 1,\ldots,n+1$ and let $A_k=\sum_{i=1}^k a_i$.
By applying the previous inequality with $a=A_{n-1}$, $b= a_n$ and $c=
a_{n+1}$, we get that
\[
\frac{A_n e^{- \alpha A_n} - A_{n-1}e^{- \alpha A_{n-1}}}{e^{-\alpha
A_{n-1}} - e^{- \alpha A_n}} > \frac{A_{n+1}e^{- \alpha A_{n+1}} - A_n
e^{- \alpha A_n}}{e^{-\alpha
A_n} - e^{- \alpha A_{n+1}}}.
\]
Recall that $\lambda_{k,m} = \sum_{i=k}^m d_{i+1}\lldots
d_m(e^{-2\alpha u_{i+1}}-e^{-2\alpha u_i})$.
The monotonicity of $\Lambda_{k,m}$ in $k$
follows easily from combining the inequality above with
the fact that if $x_1/y_1>\cdots>x_n/y_n$ and if
$y_i,c_i>0$, then
$\sum_{i=1}^{n-1}{c_i x_i}/\sum_{i=1}^{n-1}{c_i y_i} > \sum
_{i=1}^{n}{c_i x_i}/\sum_{i=1}^{n}{c_i y_i}$.
The proof of the second part of Lemma~\ref{lemreml01}
is easy and left to the reader.
The following lemma results directly from Lemma~\ref{lemreml01}.

%
\begin{lemma}\label{correml02}
With $m$ fixed and parametrization $(\nu,\alpha)$,
the quantities
$(n-1)\sum_{k=1}^m q_{k,n} (\Lambda_{k,m}-\Lambda_{m,m})^2$,
$(n-1)\sum_{k=1}^m q_{k,n} (\Lambda_{k,m}-\Lambda_{m,m})$
and the trace of
$\bolds{\Sigma}_n{(\bolds{\Sigma}_n^{-1})}_{kj}
\bolds{\Sigma}_n{(\bolds{\Sigma}_n^{-1})}_{li}$
are bounded in $O(\tilde n)$
uniformly on any compact subset of $\{T>0,\alpha>0\}$.
Therefore, \textup{(C1)(i)} and \textup{(ii$'$)} are satisfied if
$\det\mathbf{B}_n$ is of order greater than $n \tilde{n}^{1/2}$,
that is, if $\det\mathbf{B}_n^{-1}=o(n^{-1} \tilde{n}^{-1/2})$.
\end{lemma}

It is easy to see that
$ \det\mathbf{B}_n \sim2n\tilde{n}/(\nu^2v_\alpha)$
with $v_\alpha$ defined later. Indeed,
$p_{k,n}$ converges to $0$ when $k<s$,
where $s$ is the largest level $\leq m - 1$
such that $d_s$ goes to infinity and $d_{s+1},\ldots,d_{m-1}$
are fixed.
For $k=s$, $p_{s,n}$ converges to $p_{s} = 1/(d_{s+1}\lldots d_{m-1})$,
and $p_{k,n}$ converges to $p_k = (d_k - 1)/(d_k\lldots d_{m-1})$ for \mbox{$s<k<m$}.
Note that $p_s,\ldots,p_{m-1}$ are the asymptotic relative frequencies
of node ages at levels $s,\ldots,m-1$.
If $d_m$ goes to infinity, then
$v_{\alpha}=8/\sum_{k=s}^{m-1}p_k (\Lambda_k - \Lambda_m)^2$
with $\Lambda_k=\lim_n \Lambda_{k,m}$.
If $d_m$ is fixed,
\[
v_{\alpha} = 8 \Biggl(\sum_{k=s}^{m-1}p_k
(\Lambda_k - \Lambda_m)^2 - \Biggl(\sum
_{k=s}^{m-1}p_k (
\Lambda_k - \Lambda_m) \Biggr)^2\Big/d_m
\Biggr)^{-1}.
\]
Clearly, $v_\alpha>0$ because $p_{m-1}=1-1/d_{m-1}>0$ is fixed
and $\Lambda_{m-1}-\Lambda_m>0$ is easily checked.
So $\det\mathbf{B}_n$ is of order $n\tilde{n}$. The consistency and
asymptotic normality of $(\hat\nu,\hat\alpha)$
follows from applying Lemma~\ref{correml02}.

For the second part of the theorem,
we obtain the asymptotic normality of
$\sqrt{\tilde n}(\hat\gamma-\gamma,\hat\alpha-\alpha)$
through that of
$\sqrt{\tilde n}(c_1 \hat\gamma+ c_2 \hat\alpha- c_1 \gamma- c_2
\alpha)$
for every $c_1, c_2\in\mathbb R$. For this we apply the following
$\delta$-method. Its proof is similar to that of the classical
$\delta$-method [\citet{shao1999mathematical}]
and is left to the reader.

%
\begin{lemma}\label{lemreml04}
Assume that $(a_n (X_n - x),b_n (Y_n - y))^{t}$ converges
in distribution to
$N(\mathbf{0}, \bolds{\Sigma})$,
with $a_n,b_n \to\infty$, $b_n/a_n \to0$
and $\Sigma_{22} > 0$.
Suppose that $g\dvtx \mathbb{R}^2 \to\mathbb{R}$ is a continuous differentiable
function such that $\partial g/\partial y(x,y) \ne0$. Then
$b_n (g(X_n,Y_n) - g(x,y) )$ also converges
to a centered normal distribution with variance
$\Sigma_{22} ( \partial g/\partial y(x,y) )^2$.
\end{lemma}

Finally, using the classical $\delta$-method and the fact that $\sqrt
n (\hat\nu- \nu)$ is asymptotically normal, we deduce that the
asymptotic correlation between $\log\hat\gamma$ and $\log(1-\exp
^{-2\hat\alpha t_m})$ is $-1$
if $\tilde n=o(n)$.
\end{appendix}



%

\printaddresses

\end{document}